\documentclass[11pt,epsf]{article}
\usepackage{graphics,color,boxedminipage,bm,amsmath}
\usepackage{epsfig}
\usepackage{epsfig}
\usepackage{enumerate}
\usepackage{wrapfig}
\usepackage{graphics}
\usepackage{verbatim}
\usepackage{alltt}
\newtheorem{thm}{Theorem}

\newtheorem{defn}{Definition}
\newtheorem{lemma}{Lemma}

\newtheorem{conj}{Conjecture}

\newcommand{\qed}{\hspace*{.25in}\rule{.75pt}{8pt}\hspace*{.03in}\rule{2pt}{8pt}
}

\begin{document}
\large
\title{\Large \bf On degree-colorings of multigraphs\\}
\author
{\normalsize Mark K. Goldberg \\
Department of Computer Science,\\
Rensselaer Polytechnic Institute \\
Troy, NY, 12180.\\
{\tt goldbm4@rpi.edu}}

\date{December 21, 2016}
\maketitle

\begin{abstract}

A notion of degree-coloring is introduced; it captures some, but not all
properties of standard
edge-coloring. We conjecture that the smallest number of colors
needed for degree-coloring of a multigraph $G$ [the degree-coloring index
$\tau(G)$] equals $\max\{\Delta, \omega\}$, where $\Delta$ and $\omega$ are 
the maximum vertex degree in $G$ and the multigraph density, respectively. 
We prove that the conjecture holds iff $\tau(G)$ is a monotone function 
on the set of multigraphs.
\end{abstract}

\section{Introduction.}
The {\bf chromatic index} $\chi'(G)$ of a multigraph $G(V,E)$ is the minimal 
number of colors (positive integers)
that can be assigned to the edges of $G$ so that no two 
adjacent edges receive the same color.
Clearly, $ \Delta (G)\leq \chi'(G),$ where $\Delta(G)$ is the maximal vertex 
degree in $G$.
The famous result by Vizing (\cite{V-65}) establishes 
$\chi' =\chi'(G) \leq \Delta(G) + p(G)$, where $p(G)$ is the
maximal number of parallel edges in $G$. 
For graphs, in particular, $\Delta (G) \leq \chi' (G)\leq \Delta (G) +1$.
The problem of computing the exact value of the chromatic index
was proved by Holyer (\cite{Ho-81}) to be NP-hard even for cubic graphs.
It is suspected that for every multigraph with $\chi'(G) > \Delta(G) +1$,
its chromatic index is determined by the parameter $\omega(G)$, called the
multigraph {\bf density}:
$$
\omega(G)=\max_{H \subseteq G}\lceil\frac{e(H)}{\lfloor v(H)/2\rfloor}\rceil ,
$$
where $H$ is a sub-multigraph of $G$,  and $v(H)$ (resp. $e(H)$) denotes the 
number of vertices (resp. edges) in $H$.
It is easy to prove that $\omega(G) \leq \chi' (G)$ for every multigraph $G$. 
Seymour in \cite{S-79} and Stahl in \cite{Stahl-79} proved the equality
$\max (\Delta(G), \omega(G) )= \chi'^*(G)$, where $\chi'^*(G)$ is the
the fractional chromatic index of $G$, known to be  polynomially computable
(see \cite{Schein}).

The following variation of the multigraph density 
idea was considered in \cite{G-07}.
Let $\pi(F)$ denote the size of a maximum matching composed of the edges in a set $F \subseteq E$.
Denote
$\omega^*(G) = \max_{F \subseteq E} \lceil \frac{|F|}{\pi(F)}\rceil.$ 
Then, it is easy to see that
$$\omega(G) \leq \omega^*(G) \leq \chi'(G).$$
It turns out (see \cite{G-07}) that $\omega^*(G) = \max (\Delta(G), \omega(G) ).$

Conjectures connecting $\chi'(G), \omega(G)$, and $\Delta(G)$ 
were independently 
proposed by Goldberg (\cite{G-73}) and Seymour (\cite{S-79}) more than 
30 years ago 
(\cite{HaxKier},\cite{TheBook}).
(\cite{TheBook}).
Currently, the strongest variation of the conjecture (\cite{G-84}) is 
as follows: 
\begin{conj} If $\Delta(G) \not= \omega(G)$, 
then $\chi'(G) = \max (\Delta(G), \omega(G))$, 
else $\chi'(G) \leq \Delta(G) + 1 $.  
\end{conj} 

Every edge-coloring with colors $1,2,\ldots, c$ yields an assignment 
$\mu: V\rightarrow 2^{[1,c]}$, 
where  for every $x\in V$, $\mu(x)$ denotes the set 
of colors used on the edges incident to $x$.  
Given $S \subseteq V$ and $i \in [1,c]$, the set of vertices $x \in S$ such 
that $i \in \mu (x)$ is denoted $S^{(i)}(\mu).$ 
It is easy to prove that the assignment $\mu$ originated by an edge-coloring 
using colors $1,\ldots, c$ satisfies the following three conditions: 
\begin{center} 
\begin{description} 
\item[Degree condition:] $\forall x \in V(G)$, $|\mu(x)| = deg_G(x)$;
\item[Cover condition:] 	
$\forall S \subseteq V$,
$|E(S)| \leq \sum_{i=1}^c \lfloor \frac{|S^{(i)}(\mu)|}{2} \rfloor$;
\item[Matching condition:] $\forall i \in [1,c]$, the submultigraph  
induced on $V^{(i)}$ either has a perfect matching, or is empty.
\end{description}
\end{center} 
\begin{defn}
An assignment $\mu: V(G) \rightarrow 2^{[1,c]}$ satisfying the degree 
and the cover conditions is called a {\bf degree-coloring}. 
\end{defn}
Straightforward checking of the assignment presented in the Figure below 
shows that the assignment is a degree-coloring of the multicycle $C$. 
However, it is not originated by any edge-coloring of $C$, since the 
submultigraph of $C$ induced on $V^{(6)}$ has no perfect matching.

\begin{center}
\psfig{figure=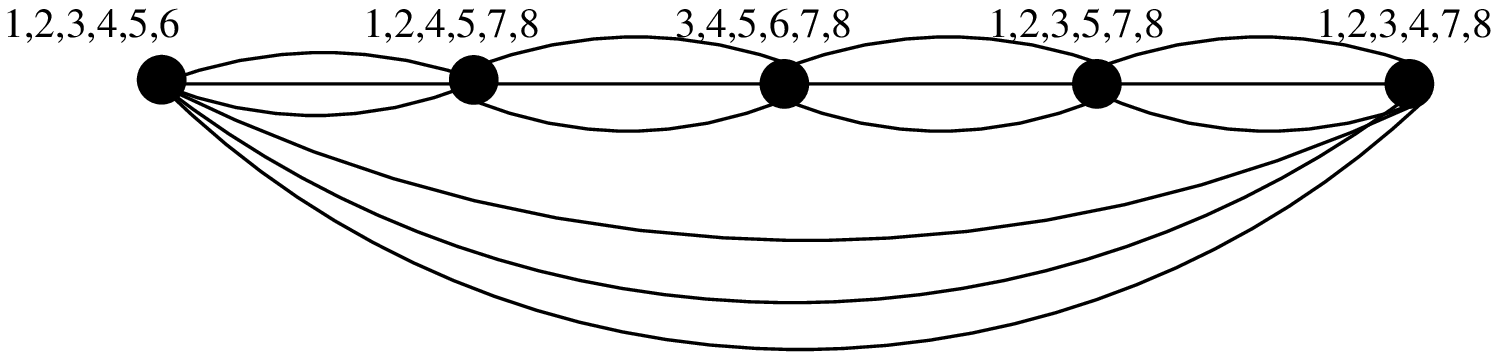,height=0.8in}
\end{center}
Let $\tau (G)$ denote the smallest integer $c$ 
for which a degree-coloring of $G$ exists.
It is easy to prove 
\begin{lemma} \label{lower-bound}
$\max(\Delta(G), \omega(G)) \leq \tau(G) \leq \chi'(G)$.
\end{lemma}

\begin{conj}\label{pre-conjecture} (the $\tau$-conjecture):
For every multigraph $G$, $\tau(G) = \max(\Delta(G), \omega(G))$.
\end{conj}
A  real-valued function $\kappa (G)$ defined on  the set of
multigraphs is called {\bf monotone} if for any multigraph $G$ and
any submultigraph $H \subseteq G$, $\kappa(H) \leq \kappa(G)$.
Clearly, $\Delta(G)$ and $\omega(G)$ are monotone functions.

\begin{conj} \label{monotone-conjecture}
The degree-coloring index $\tau(G)$ is a monotone function on multigraphs.
\end{conj}

It is easy to see that Conjecture \ref{pre-conjecture} 
implies Conjecture \ref{monotone-conjecture}.
We prove in this paper that the reverse is also true: the monotonicity 
of $\tau(G)$ implies conjecture \ref{pre-conjecture}.

We use the standard graph-theoretical terminology which can be 
found in \cite{We}.

\section{Monotonicity of $\tau(G)$ and the $\tau$-conjecture.}
It is easy to construct a $\tau(G)$-degree-coloring for a regular 
multigraph $G$ with $\omega(G) \leq \Delta(G)$.
\begin{lemma}\label{upper-bound}
If $G$ is a $\Delta$-regular multigraph, and 
$\omega(G)\leq\Delta$,  then $\tau(G) = \Delta.$
\end{lemma}
{\bf Proof.} From the definition, $\tau(G) \geq \Delta$. 
Consider the following assignment:
$$
\forall x \in V(G), ~
\mu(x) = \{1, 2, \ldots ,\Delta\}.
$$
Given $S \subseteq V(G)$, 
$\forall i \in [1, \Delta]$, $S^{(i)}(\mu) = S$. Thus,
$$\sum_{i=1}^{\Delta}\lfloor\frac{|S^{(i)}(\mu)|}{2} \rfloor 
= \lfloor \frac{|S|}{2}\rfloor\Delta .$$
Since $\omega(G) \leq \Delta$, for any $S \subseteq V$, 
$\lfloor \frac{|S|}{2}\rfloor \Delta \geq 
\lfloor \frac{|S|}{2}\rfloor \omega (G) \geq |E(S)|$ implying 
$\tau( G) = \Delta.$ \qed

Constructing a degree-coloring for a non-regular multigraph can be done via 
operation Regularization which, for every multigraph $G$, 
creates a regular multigraph $R(G)$ containing $G$ as 
an induced sub-multigraph.

{\bf Regularization:} If a multigraph $G$ is regular and 
$\omega(G)\leq \Delta(G)$, then $R(G) = G$; else
\begin{enumerate} 
\item generate a disjoint isomorphic copy $G' = (V', E')$  of $G(V,E)$ 
with an isomorphic mapping $f: V\rightarrow V'$ from $G$ onto $G'$; 
\item 
let $V(R(G)) = V \cup V'$ and initialize $E(R(G))$ by setting 
$E(R(G)) = E(G) \cup E(G')$;
\item
$\forall x \in V$, add $\max(\Delta(G),\omega(G))-deg_G(x)$ new edges 
$xf(x)$ to $E(R(G))$.
\end{enumerate}
\begin{lemma} $\forall G$, 
$\omega(G) \leq \omega(R(G)) \leq \max(\omega(G), \Delta(G))$ and 
$\Delta(R(G)) = \max (\Delta(G), \omega(G)).$
\end{lemma}
{\bf Proof.} If $G = R(G)$, the lemma is obvious. Let $G\not= R(G)$. 
Denote $\Delta = \Delta(G),$ $\omega = \omega(G)$, and 
$\rho = \max(\Delta, \omega )$.
Obviously, $\Delta(R(G)) = \rho$
and $\omega \leq \omega(R(G))$.  

To prove $\omega(R(G)) \leq \rho$, denote $R= R(G)$, $V(R) = V_1 \cup V_2$, 
where $V_1 = V(G)$ and
$V_2 = V(G')$. Let $f$ be an isomorphic mapping from $V_1$ onto $V_2$.
Given $S \subseteq V(R)$, let
$S_1 = S \cap V_1$, $S_2 =  S\cap V_2$,
$S' = S_1 \cap f^{-1} (S_2)$,  and $S'' =  S_2 \cap f(S_1)$.
Note that $|S'| = |S''|$ and $|E(S')| = |E(S'')|.$
\begin{center}
\psfig{figure=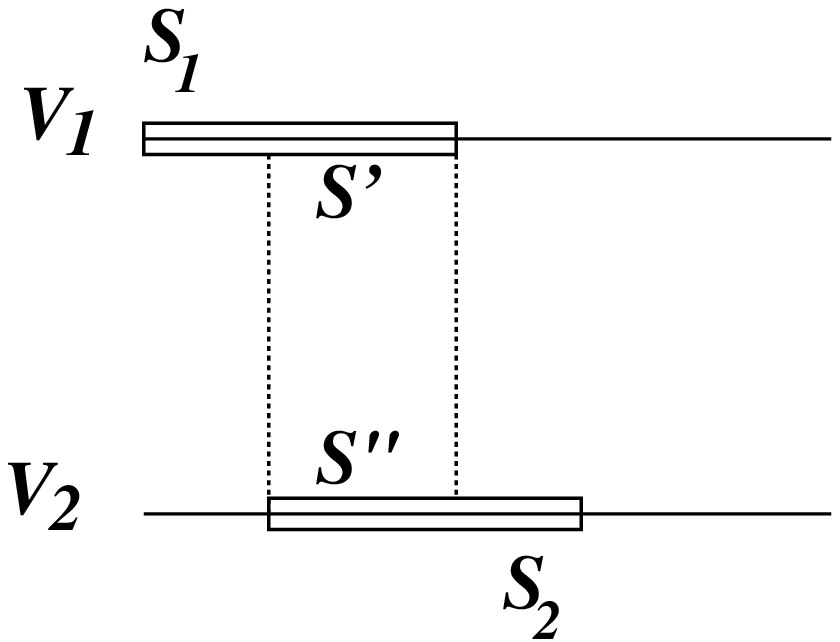,height=1.8in}
\end{center}

\begin{eqnarray*}
\mbox{Then~} |E(S) |&=& |E(S_1)| + |E(S_2)| + \sum_{x\in S'} (\rho - deg_G(x))\\
&=& |E(S_1-S')| + |E(S_1 -S', S')| + |E(S')| + \\
&& |E(S_2-S'')|+|E(S_2 -S'', S'')| + |E(S'')|  + 
|S'| \rho - \sum_{x\in S'} deg_G(x).
\end{eqnarray*}
It is easy to check that
\begin{eqnarray*}
|E(S_1 -S', S')| + |E(S')|+ |E(S_2 -S'', S'')| + |E(S'')| =\\
|E(S_1 -S', S')| + |E(S_2 -S'', S'')| + 2 |E(S')| \leq 
\sum_{x\in S'} deg_G(x)),
\end{eqnarray*}
which yields the following upper bound
\begin{eqnarray*}
|E(S)| &\leq &|E(S_1 -S')|+|E(S_2 -S'')| + |S'| \rho \\ 
&\leq &
\lfloor \frac{|S_1| -|S'|}{2} \rfloor \rho + 
\lfloor \frac{|S_2| -|S''|}{2} \rfloor \rho + |S'|\rho .
\end{eqnarray*}
To prove
$$
\lfloor \frac{|S_1| -|S'|}{2} \rfloor \rho + 
\lfloor \frac{|S_2| -|S''|}{2} \rfloor \rho + |S'|\rho  
\leq \lfloor\frac{|S_1|+|S_2|}{2} \rfloor \rho ,$$
note that it is straightforward if $|S_1| + |S_2|$ is even.
If $|S_1| + |S_2|$ is odd, one out of two integers
$|S_1| -|S'|$
and $|S_2| -|S''|$ is even and one is odd.
Thus, 
$$
\lfloor \frac{|S_1| -|S'|}{2} \rfloor \rho + 
\lfloor \frac{|S_2| -|S''|}{2} \rfloor \rho + |S'| \rho =
\frac{|S_1|}{2} \rho + \frac{|S_2|}{2}  \rho - 
\frac{1}{2} \rho.$$
Since $|S_1| + |S_2|$ is odd,
$$
\lfloor\frac{|S_1|+|S_2|}{2} \rfloor \rho = 
\frac{|S_1|+|S_2|}{2} \rho - \frac{1}{2} \rho ,
$$
which implies the result. \qed

\begin{thm} If function $\tau(G)$ is monotone on the set of all multigraphs,
then for any multigraph $G$, 
$$ \tau(G) = \max\{\Delta(G), \omega(G)\}.$$
\end{thm}
{\bf Proof.} By Lemma \ref{lower-bound},
$  \max\{\Delta(G), \omega(G)\}\leq \tau(G) $. On the other hand, since 
$ G \subseteq R(G)$,  it follows from Lemma \ref{upper-bound} that
$ \tau(G) \leq \tau(R(G)) = \max\{\Delta(G), \omega(G)\}.$
\qed

\end{document}